
\input amstex
\input amsppt.sty
\hsize 30pc
\vsize 47pc
\def\nmb#1#2{#2}         
\def\cit#1#2{\ifx#1!\cite{#2}\else#2\fi}
\def\totoc{}             
\def\idx{}               
\def\ign#1{}             
\redefine\o{\circ}
\define\X{\frak X}
\define\al{\alpha}

\define\ep{\varepsilon}
\define\ze{\zeta}

\define\rh{\rho}
\define\si{\sigma}

\define\ph{\varphi}

\define\om{\omega}

\define\Si{\Sigma}

\define\Om{\Omega}
\redefine\i{^{-1}}
\define\row#1#2#3{#1_{#2},\ldots,#1_{#3}}
\define\x{\times}

\define\Ad{\operatorname{Ad}}

\redefine\L{{\Cal L}}

\define\g{{\frak g}}

\def\today{\ifcase\month\or
 January\or February\or March\or April\or May\or June\or
 July\or August\or September\or October\or November\or December\fi
 \space\number\day, \number\year}
\topmatter
\title  Basic differential forms for actions of Lie groups
\endtitle
\author
Peter W\. Michor  \endauthor
\affil
Erwin Schr\"odinger International Institute of Mathematical Physics,
Wien, Austria \\
Institut f\"ur Mathematik, Universit\"at Wien, Austria
\endaffil
\address
P\. W\. Michor: Institut f\"ur Mathematik, Universit\"at Wien,
Strudlhofgasse 4, A-1090 Wien, Austria
\endaddress
\email MICHOR\@ESI.AC.AT \endemail
\date {May 2, 1994} \enddate
\thanks Supported by Project P 10037--PHY
of `Fonds zur F\"orderung der wissenschaftlichen Forschung'.
\endthanks
\subjclass Orbits, sections, basic differential forms\endsubjclass
\abstract A section of a Riemannian $G$-manifold $M$ is a closed
submanifold $\Si$ which meets each orbit orthogonally.
It is shown
that the algebra of $G$-invariant differential forms on $M$ which are
horizontal in the sense that they kill every vector which is tangent
to some orbit, is isomorphic to the algebra of those differential forms on
$\Si$ which are invariant with respect to the generalized Weyl group
of this orbit, under some condition.
\endabstract
\endtopmatter

\document

\heading Table of contents \endheading

\noindent 1. Introduction \leaders \hbox to
1em{\hss .\hss }\hfill {\eightrm 1}\par
\noindent 2. Basic differential forms \leaders \hbox to
1em{\hss .\hss }\hfill {\eightrm 2}\par
\noindent 3. Representations \leaders \hbox to
1em{\hss .\hss }\hfill {\eightrm 3}\par
\noindent 4. Proof of the main theorem \leaders \hbox to
1em{\hss .\hss }\hfill {\eightrm 7}\par
\noindent 5. Basic versus equivariant cohomology \leaders \hbox to
1em{\hss .\hss }\hfill {\eightrm 8}\par

\head\totoc\nmb0{1}. Introduction \endhead
A section of a Riemannian $G$-manifold $M$ is a closed
submanifold $\Si$ which meets each orbit orthogonally. This notion
was introduced by Szenthe \cit!{22}, \cit!{23}, in slightly different
form by Palais and Terng in \cit!{17}, \cit!{18}.
The case of linear representations was
considered by Bott and Samelson \cit!{4}, Conlon \cit!{9},
and then by Dadok \cit!{10} who called
representations admitting sections polar representations and
completely classified all polar representations of connected Lie
groups. Conlon \cit!{8} considered Riemannian manifolds admitting
flat sections. We follow here the notion of Palais and Terng.

If $M$ is a Riemannian $G$-manifold which admits a section $\Si$
then the trace on $\Si$ of the $G$-action is a discrete group action
by the generalized Weyl group $W(\Si)=N_G(\Si)/Z_G(\Si)$. Palais and
Terng \cit!{17} showed that then the algebras of invariant smooth
functions coincide
$C^\infty(M,\Bbb R)^G\cong C^\infty(\Si,\Bbb R)^{W(\Si)}$.

In this paper we will extend this result to the algebras of
differential forms. Our aim is to show that pullback along the
embedding $\Si\to M$ induces an isomorphism
$\Om_{\text{hor}}^p(M)^G\cong \Om^p(\Si)^{W(\Si)}$ for each $p$,
where a differential form $\om$ on $M$ is called \idx{\it horizontal}
if it kills each vector tangent to some orbit. For each point $x$ in
$M$, the slice representation of the isotropy group $G_x$ on the
normal space $T_x(G.x)^\bot$ to the tangent space to the orbit
through $x$ is a polar representation. The first step is to show that
the result holds for polar representations. This is done in theorem
\nmb!{3.6} for polar representations whose generalized Weyl group is
really a Coxeter group, is generated by reflections. Every polar
representation of a connected Lie group has this property. The method used
there is inspired by Solomon \cit!{21}. Then the general result is
proved under the assumption that each slice representation has a
Coxeter group as a generalized Weyl group.

I want to thank D\. Alekseevsky for introducing me to the beautiful
results of Palais and Terng, and him and A\. Onishchik for many
discussions about this and related topics.

\head\totoc\nmb0{2}. Basic differential forms \endhead

\subhead\nmb.{2.1}. Basic differential forms \endsubhead
Let $G$ be a Lie group with Lie algebra $\g$, multiplication
$\mu:G\x G\to G$, and for $g\in G$ let $\mu_g, \mu^g:G\to G$ denote
the left and right translation.

Let $\ell:G\x M\to M$ be a left action of the Lie group $G$ on a smooth
manifold $M$. We consider the partial mappings $\ell_g:M\to M$ for
$g\in G$ and $\ell^x:G\to M$ for $x\in M$ and the fundamental vector
field mapping $\ze:\g\to \X(M)$ given by $\ze_X(x)=T_e(\ell^x)X$.
Since $\ell$ is a left action, the negative $-\ze$ is a Lie
algebra homomorphism.

A differential form $\ph\in \Om^p(M)$ is called
\idx{\it $G$-invariant} if $(\ell_g)^*\ph=\ph$ for all $g\in G$ and
\idx{\it horizontal} if $\ph$ kills each vector tangent to a
$G$-orbit: $i_{\ze_X}\ph=0$ for all $X\in\g$.
We denote by $\Om^p_{\text{hor}}(M)^G$ the space of all horizontal
$G$-invariant $p$-forms on $M$. They are also called \idx{\it basic
forms}.

\proclaim{\nmb.{2.2}. Lemma} Under the exterior differential
$\Om_{\text{hor}}(M)^G$ is a subcomplex of $\Om(M)$.
\endproclaim

\demo{Proof}
If $\ph\in\Om_{\text{hor}}(M)^G$ then the exterior derivative
$d\ph$ is clearly $G$-invariant. For $X\in\g$ we have
$$i_{\ze_X}d\ph = i_{\ze_X}d\ph + di_{\ze_X}\ph = \L_{\ze_X}\ph=0,$$
so $d\ph$ is also horizontal.
\qed\enddemo

\subhead\nmb.{2.3}. Sections \endsubhead
Let $M$ be a connected complete Riemannian manifold and let $G$ be a
Lie group which acts isometrically on $M$ from the left.
A connected closed smooth
submanifold $\Si$ of $M$ is called a \idx{\it section} for the
$G$-action, if it meets all $G$-orbits orthogonally.

Equivalently we require that $G.\Si=M$ and that for each $x\in\Si$
and $X\in \g$ the fundamental vector field $\ze_X(x)$ is orthogonal
to $T_x\Si$.

We only remark here
that each section is a totally geodesic submanifold and is given by
$\exp(T_x(x.G)^\bot)$ if $x$ lies in a principal orbit.

If we put $N_G(\Si):= \{g\in G: g.\Si=\Si\}$ and $Z_G(\Si):=
\{g\in G: g.s=s \text{ fo all }s\in \Si\}$, then the quotient
$W(\Si):= N_G(\Si)/Z_G(\Si)$ turns out to be a discrete group acting
properly on $\Si$. It is called the generalized Weyl group of the
section $\Si$.

See \cit!{17} or \cit!{18}
for more information on sections and their generalized
Weyl groups.

\proclaim{\nmb.{2.4}. Main Theorem} Let $M\x G\to M$ be a proper
isometric right action of a Lie group $G$ on a a smooth Riemannian
manifold $M$, which admits a section $\Si$. Let us assume that
\roster
\item For each $x\in\Si$ the slice representation
     $G_x\to O(T_x(G.x)^\bot)$ has a generalized Weyl group which
     is a reflection group (see section \nmb!{3}).
\endroster

Then the restriction of
differential forms induces an isomorphism
$$\Om^p_{\text{hor}}(M)^G @>{\cong}>> \Om^p(\Si)^{W(\Si)}$$
between the space of horizontal $G$-invariant differential forms on
$M$ and the space of all differential forms on $\Si$ which are
invariant under the action of the generalized Weyl group $W(\Si)$ of
the section $\Si$.
\endproclaim

The proof of this theorem will take up the rest of this paper.
According to Dadok \cit!{10}, remark after Proposition 6, for any
polar representation of a connected Lie group
the generalized Weyl group $W(\Si)$ is a
reflection group, so condition \thetag1 holds if we assume that:
\roster
\item [2] Each isotropy group $G_x$ is connected.
\endroster

\head\totoc\nmb0{3}. Representations \endhead

\subhead\nmb.{3.1}. Invariant functions \endsubhead
Let $G$ be a reductive Lie group and let $\rh:G\to GL(V)$ be a
representation in a finite dimensional real vector space
$V$.

According to a classical theorem of Hilbert (as extended by Nagata
\cit!{13}, \cit!{14}), the algebra of
$G$-invariant polynomials $\Bbb R[V]^G$ on $V$ is finitely generated
(in fact finitely presented),
so there are $G$-invariant homogeneous polynomials $f_1,\dots,f_m$ on
$V$ such that each invariant polynomial $h\in \Bbb R[V]^G$ is of the
form $h= q(f_1,\dots,f_m)$ for a polynomial $q\in\Bbb R[\Bbb R^m]$.
Let $f=(f_1,\dots,f_m):V\to \Bbb R^m$, then this means that the
pullback homomorphism $f^*:\Bbb R[\Bbb R^m]\to \Bbb R[V]^G$ is
surjective.

D\. Luna proved in \cit!{12}, that the pullback homomorphism
$f^*:C^\infty(\Bbb R^m,\Bbb R)\to C^\infty(V,\Bbb R)^G$ is also surjective
onto the space of all smooth functions on $V$ which are constant on
the fibers of $f$. Note that the polynomial mapping $f$ in this case
may not separate the $G$-orbits.

G\. Schwarz proved already in \cit!{20}, that if $G$ is a compact Lie group
then the pullback homomorphism
$f^*:C^\infty(\Bbb R^m,\Bbb R)\to C^\infty(V,\Bbb R)^G$ is actually
surjective onto the space of $G$-invariant smooth functions.
This result implies in particular that $f$ separates the $G$-orbits.

\proclaim{\nmb.{3.2}. Lemma} Let $\ell\in V^*$ be a linear functional
on a finite dimensional vector space $V$, and let
$f\in C^\infty(V,\Bbb R)$ be a smooth function which vanishes on the
kernel of $\ell$, so that $f|\ell\i(0)=0$. Then there is a unique
smooth function $g$ such that $f=\ell.g$
\endproclaim

\demo{Proof}
Choose coordinates $x^1,\dots,x^n$ on $V$ with $\ell=x^1$. Then
$f(0,x^2,\dots,x^n)=0$ and we have
$f(x^1,\dots,x^n)=\int_0^1\partial_1
f(tx^1,x^2,\dots,x^n)dt.x^1=g(x^1,\dots,x^n).x^1.$
\qed\enddemo

\proclaim{\nmb.{3.3}. Question} Let $\rh:G\to GL(V)$ be a
representation of a compact Lie group in a finite dimensional vector
space $V$. Let $f=(f_1,\dots,f_m):V\to \Bbb R^m$ be the polynomial
mapping whose components $f_i$ are a minimal set of
homogeneous generators for the
algebra $\Bbb R[V]^G$ of invariant polynomials.

We consider the pullback homomorphism
$f^*:\Om^p(\Bbb R^m)\to \Om^p(V)$. Is it surjective onto the space
$\Om^p_{\text{hor}}(V)^G$ of $G$-invariant
horizontal smooth $p$-forms on $V$?
\endproclaim

See remark \nmb!{3.7} for a class of representations where the
answer is yes.

In general the answer is no. A counterexample is the following:
Let the cyclic group $\Bbb Z_n=\Bbb Z/n\Bbb Z$ of order $n$, viewed
as the group of $n$-th roots of unity, act on $\Bbb C=\Bbb R^2$ by
complex multiplication. A generating system of
polynomials consists of $f_1=|z|^2$, $f_2=\operatorname{Re}(z^n)$,
$f_3=\operatorname{Im}(z^n)$. But then each $df_i$ vanishes at 0 and
there is no chance to have the horizontal
invariant volume form $dx\wedge dy$ in
$f^*\Om(\Bbb R^3)$.

\subhead\nmb.{3.4}. Polar representations \endsubhead
Let $G$ be a compact Lie group and let $\rh:G\to GL(V)$ be an
orthogonal representation in a finite dimensional real vector space
$V$ which admits a section $\Si$. Then the section turns out to be a
linear subspace and the representation is called a \idx{\it polar
representation}, following Dadok \cit!{10}, who gave a complete
classification of all polar representations of connected Lie groups.
They were called variationally complete representations by Conlon
\cit!{9} before.

\proclaim{\nmb.{3.5}. Theorem} {\rm \cit!{17}, 4.12, or \cit!{24},
theorem D.}
Let $\rh:G\to GL(V)$ be a polar representation of
a compact Lie group $G$, with section $\Si$ and generalized Weyl
group $W=W(\Si)$.

Then the algebra $\Bbb R[V]^G$ of $G$-invariant polynomials on $V$ is
isomorphic to the algebra $\Bbb R[\Si]^W$ of $W$-invariant
polynomials on the section $\Si$, via the restriction mapping
$f\mapsto f|\Si$.
\endproclaim

\proclaim{\nmb.{3.6}. Theorem}
Let $\rh:G\to GL(V)$ be a polar representation of
a compact Lie group $G$, with section $\Si$ and generalized Weyl
group $W=W(\Si)$. Let us suppose that $W=W(\Si)$ is generated by
reflections (a reflection group or Coxeter group).

Then the pullback to $\Si$ of differential forms induces an isomorphism
$$\Om^p_{\text{hor}}(V)^G @>{\cong}>> \Om^p(\Si)^{W(\Si)}.$$
\endproclaim

According to Dadok \cit!{10}, remark after proposition 6, for any
polar representation the generalized Weyl group $W(\Si)$ is a
reflection group. This theorem is true for polynomial
differential forms, and also for real analytic differential forms, by
essentially the same proof.

\demo{Proof}
Let $f_1,\dots,f_n$ be a minimal set of homogeneous generators of the
algebra $\Bbb R[\Si]^W$ of $W$-invariant polynomials on $\Si$. Then
this is a set of algebraically independent polynomials, $n=\dim \Si$,
and their degrees $d_1,\dots,d_n$ are uniquely determined up to
order. We even have (see \cit!{11})
\roster
\item $d_1\dots d_n=|W|$, the order of $W$,
\item $d_1+\dots +d_n=n+N$, where $N$ is the number of reflections in
       $W$,
\item $\prod_{i=1}^n(1+(d_i-1)t)=a_0+a_1t +\dots+a_n t^n$, where
       $a_i$ is the number of elements in $W$ whose fixed point set
       has dimension $n-i$.
\endroster
Let us consider the mapping $f=(f_1,\dots,f_n):\Si\to \Bbb R^n$ and
its Jacobian $J(x)=\det(df(x))$. Let $x^1,\dots,x^n$ be coordinate
functions in $\Si$. Then for each $\si\in W$ we have
$$\align
J.dx^1\wedge \dots \wedge dx^n &= df_1\wedge \dots \wedge df_n
     = \si^*(df_1\wedge \dots \wedge df_n) \\
& = (J\o\si)\si^* (dx^1\wedge \dots \wedge dx^n)
     = (J\o\si)\det(\si)(dx^1\wedge \dots \wedge dx^n),\\
J\o\si &= \det(\si\i) J.\tag 4
\endalign$$
If $J(x)\ne 0$, then in a neighborhood of $x$ the mapping $f$ is a
diffeomorphism by the inverse function theorem, so that the 1-forms
$df_1,\dots,df_n$ are a local coframe there.
Since the generators $f_1,\dots,f_n$ are algebraically independent
over $\Bbb R$, $J\ne 0$.
Since $J$ is a polynomial of degree $(d_1-1)+\dots+(d_n-1)=N$ (see
\therosteritem2), the set $U=\Si\setminus J\i(0)$ is open and dense
in $\Si$, and $df_1,\dots,df_n$ form a coframe on $U$.

Now let $(\si_\al)_{\al=1,\dots,N}$ be the set of reflections in $W$, with
reflection hyperplanes $H_\al$. Let $\ell_\al\in \Si^*$ be a linear
functional with $H_\al=\ell\i(0)$. If $x\in H_\al$ we have
$J(x)=\det(\si_\al)J(\si_\al.x)=-J(x)$, so that $J|H_\al=0$ for each
$\al$, and by
lemma \nmb!{3.2} we have
$$J=c.\ell_1\dots\ell_N.\tag 5$$
Since $J$ is a polynomial of degree $N$, $c$ must be a constant.
Repeating the last argument for an arbitrary function $g$ and using
\thetag5, we get:
\roster
\item [6] If $g\in C^\infty(\Si,\Bbb R)$ satisfies
     $g\o \si=\det(\si\i)g$ for each $\si\in W$, we have
     $g=J.h$ for $h\in C^\infty(\Si,\Bbb R)^W$.
\item {\bf Claim.} Let $\om\in\Om^p(\Si)^W$. Then we have
$$
\om = \sum_{j_1<\dots<j_p}\om_{j_1\dots j_p}
     df_{j_1}\wedge \dots \wedge df_{j_p},
$$
     where $\om_{j_1\dots j_p} \in C^\infty(\Si,\Bbb R)^W$.
\endroster
Since $df_1,\dots,df_n$ form a coframe on the $W$-invariant dense open set
$U=\{x:J(x)\ne0\}$, we have
$$
\om|U = \sum_{j_1<\dots<j_p}g_{j_1\dots j_p}
     df_{j_1}|U\wedge \dots \wedge df_{j_p}|U
$$
for $g_{j_1\dots j_p} \in C^\infty(U,\Bbb R)$. Since $\om$ and all
$df_i$ are $W$-invariant we may replace $g_{j_1\dots j_p}$ by
$$
\tfrac1{|W|}\sum_{\si\in W}
     g_{j_1\dots j_p}\o \si\in C^\infty(U,\Bbb R)^W,
$$
or assume without loss that
$g_{j_1\dots j_p}\in C^\infty(U,\Bbb R)^W$.

Let us choose now a form index $i_1<\dots<i_p$ with
$\{i_{p+1}<\dots<i_n\}=\{1,\dots,n\}\setminus\{i_1<\dots<i_p\}$. Then
for some sign $\ep=\pm 1$ we have
$$\align
\om|U\wedge df_{i_{p+1}}\wedge \dots\wedge df_{i_n} &= \ep.
     g_{i_1\dots i_p}. df_1\wedge \dots\wedge df_n \\
&= \ep.  g_{i_1\dots i_p}.J. dx^1\wedge \dots\wedge dx^n, \text{ and} \\
\om \wedge df_{i_{p+1}}\wedge \dots\wedge df_{i_n} &=
     \ep.k_{i_1\dots i_p} dx^1\wedge \dots\wedge dx^n\tag8
\endalign$$
for a function $k_{i_1\dots i_p}\in C^\infty(\Si,\Bbb R)$. Thus
$$
k_{i_1\dots i_p}|U = g_{i_1\dots i_p}.J|U.\tag9
$$
Since $\om$ and each $df_i$ is $W$-invariant, from \thetag8 we get
$k_{i_1\dots i_p}\o \si = \det(\si\i)k_{i_1\dots i_p}$ for each
$\si\in W$. But then by \thetag6 we have
$k_{i_1\dots i_p}=\om_{i_1\dots i_p}.J$ for unique
$\om_{i_1\dots i_p}\in C^\infty(\Si,\Bbb R)^W$, and \thetag9 then
implies $\om_{i_1\dots i_p}|U=g_{i_1\dots i_p}$, so that the claim
\thetag7 follows since $U$ is dense.

Now we may finish the proof of the theorem. Let $i:\Si\to V$ be the
embedding. By theorem \nmb!{3.5} the algebra
$\Bbb R[V]^G$ of $G$-invariant polynomials on $V$ is
isomorphic to the algebra $\Bbb R[\Si]^W$ of $W$-invariant
polynomials on the section $\Si$, via the restriction mapping $i^*$.
Choose polynomials $\tilde f_1,\dots\tilde f_n\in \Bbb R[V]^G$ with
$\tilde f_i\o i = f_i$ for all $i$. Put
$\tilde f=(\tilde f_1,\dots,\tilde f_n):V\to \Bbb R^n$.
In the setting of claim \thetag7,
use the theorem of G. Schwarz (see \nmb!{3.1}) to find
$h_{i_1,\dots,i_p}\in C^\infty(\Bbb R^n,\Bbb R)$ with
$h_{i_1,\dots,i_p}\o f = \om_{i_1,\dots,i_p}$ and consider
$$\tilde\om = \sum_{j_1<\dots<j_p}(h_{j_1\dots j_p}\o \tilde f)
     d\tilde f_{j_1}\wedge \dots \wedge d\tilde f_{j_p},$$
which is in $\Om^p_{\text{hor}}(V)^G$ and satifies
$i^*\tilde\om=\om$.

Thus the mapping
$i^*:\Om^p_{\text{hor}}(V)^G\to \Om^p_{\text{hor}}(\Si)^W$
is surjective. It is also injective:
Let $\om\in \Om^p_{\text{hor}}(V)^G$ with $i^*\om=0$. Then
for a regular point $x\in \Si$ we have $\om_x=0$ since it vanishes on
vectors orthogonal to the orbit $G.x$ by $(i^*\om)_x=0$, and on
vectors tangential to the orbit $G.x$ by horizontality. By
$G$-invariance then $\om$ vanishes along the whole orbit $G.x$. Since
regular orbits are dense in $V$, $\om=0$.
\qed\enddemo

\subhead\nmb.{3.7}. Remark \endsubhead
The proof of theorem \nmb!{3.6} shows that the answer to question
\nmb!{3.3} is yes for the representations treated in \nmb!{3.6}.

\proclaim{\nmb.{3.8}. Corollary}
Let $\rh:G\to O(V,\langle\quad,\quad\rangle)$ be an orthogonal polar
representation of a compact Lie group $G$, with section $\Si$ and
generalized Weyl group $W=W(\Si)$. Let us suppose that $W=W(\Si)$ is
generated by reflections (a reflection group or Coxeter group). Let
$B\subset V$ be an open ball centered at 0.

Then the restriction of differential forms induces an isomorphism
$$\Om^p_{\text{hor}}(B)^G @>{\cong}>> \Om^p(\Si\cap B)^{W(\Si)}.$$
\endproclaim

\demo{Proof} Check the proof of \nmb!{3.6} or use the following
argument.
Suppose that $B=\{v\in V:|v|< 1\}$ and consider a smooth
diffeomorphism $f:[0,1)\to [0,\infty)$ with $f(t)=t$ near 0. Then
$g(v):= \frac {f(|v|)}{|v|}v$ is a $G$-equivariant diffeomorphism
$B\to V$ and by \nmb!{3.6} we get:
$$
\Om^p_{\text{hor}}(B)^G @>{(g\i)^*}>> \Om^p_{\text{hor}}(V)^G
     @>{\cong}>> \Om^p(\Si)^{W(\Si)} @>{g^*}>>
     \Om^p(\Si\cap B)^{W(\Si)}.
\qed$$
\enddemo

\head\totoc\nmb0{4}. Proof of the main theorem \endhead

Let us assume that we are in the situation of the main theorem
\nmb!{2.4}, for the rest of this section.

\subhead\nmb.{4.1}  \endsubhead
Let $i:\Si\to M$ be the embedding of the section.
It clearly induces a linear mapping
$i^*:\Om^p_{\text{hor}}(M)^G \to \Om^p(\Si)^{W(\Si)}$
which is injective by the following argument: Let
$\om\in\Om^p_{\text{hor}}(M)^G$ with $i^*\om=0$.
For $x\in \Si$ we have $i_X\om_x=0$ for $X\in T_x\Si$ since
$i^*\om=0$, and also for $X\in T_x(G.x)$ since $\om$ is horizontal.
Let $x\in\Si\cap M_{\text{reg}}$ be a regular point, then $T_x\Si =
(T_x(G.x))^\bot$ and so $\om_x=0$. This holds along the whole orbit
through $x$ since $\om$ is $G$-invariant. Thus
$\om|M_{\text{reg}}=0$, and since $M_{\text{reg}}$ is dense in $M$,
$\om=0$.

So it remains to show that $i^*$ is surjective.

\subhead\nmb.{4.2}  \endsubhead
For $x\in M$ let $S_x$ be a (normal) slice and
$G_x$ the isotropy group, which acts on the slice. Then $G.S_x$ is
open in $M$ and $G$-equivariantly diffeomorphic to the associated
bundle $G\to G/G_x$ via
$$\CD
G\x S_x @>{q}>> G\x_{G_x}S_x @>{\cong}>> G.S_x \\
@.               @VVV                    @VVrV \\
@.              G/G_x @>{\cong}>>        G.x,
\endCD$$
where $r$ is the projection of a tubular neighborhood.
Since $q:G\x S_x \to G\x_{G_x}S_x$ is a principal $G_x$-bundle with
principal right action $(g,s).h=(gh,h\i.s)$, we have an isomorphism
$q^*:\Om(G\x_{G_x}S_x) \to \Om_{G_x-\text{hor}}(G\x S_x)^{G_x}$.
Since $q$ is also $G$-equivariant for the left $G$-actions, the
isomorphism $q^*$ maps the subalgebra
$\Om_{\text{hor}}^p(G.S_x)^G\cong \Om_{\text{hor}}^p(G\x_{G_x}S_x)^G$
of $\Om(G\x_{G_x}S_x)$ to the subalgebra
$\Om_{G_x-\text{hor}}^p(S_x)^{G_x}$ of
$\Om_{G_x-\text{hor}}(G\x S_x)^{G_x}$. So we have proved:

\proclaim{Lemma} In this situation
there is a canonical isomorphism
$$\Om_{\text{hor}}^p(G.S_x)^G @>{\cong}>>
     \Om_{G_x-\text{hor}}^p(S_x)^{G_x}$$
which is given by pullback along the embedding $S_x\to G.S_x$.
\endproclaim

\subhead\nmb.{4.3}. Rest of the proof of theorem \nmb!{3.6} \endsubhead
Now let us consider $\om\in\Om^p(\Si)^{W(\Si)}$. We want to construct
a form $\tilde\om\in\Om^p_{\text{hor}}(M)^G$ with
$i^*\tilde\om =\om$. This will finish the proof of theorem
\nmb!{3.6}.

Choose $x\in \Si$ and an open ball $B_x$ with center $0$ in $T_xM$
such that the Riemannian exponential mapping $\exp_x:T_xM\to M$ is a
diffeomorphism on $B_x$. We considernow the compact isotropy group
$G_x$ and the slice representation
$\rh_x:G_x\to O(V_x)$, where
$V_x=\operatorname{Nor}_x(G.x)=(T_x(G.x))^\bot\subset T_xM$
is the normal space to the orbit. This is a polar representation with
sextion $T_x\Si$, and its generalized Weyl group is given by
$W(T_x\Si)\cong N_G(\Si)\cap G_x/Z_G(\Si)= W(\Si)_x$ (see
\cit!{17}\ign{9.4}) and it is a Coxeter
group by assumption \therosteritem1 in \nmb!{3.6}.
Then $\exp_x: B_x\cap V_x\to S_x$ is a diffeomorphism onto a slice
and $\exp_x: B_x\cap T_x\Si \to \Si_x\subset\Si$ is a diffeomorphism
onto an open neighborhood $\Si_x$ of $x$ in the section $\Si$.

Let us now consider the pullback
$(\exp|B_x\cap T_x\Si)^*\om \in \Om^p(B_x\cap T_x\Si)^{W(T_x\Si)}$.
By corollary \nmb!{3.8} there exists a unique form
$\ph^x\in\Om^p_{G_x-\text{hor}}(B_x\cap V_x)^{G_x}$ such that
$i^*\ph^x=(\exp|B_x\cap T_x\Si)^*\om$, where $i_x$ is the embedding.
Then we have
$$((\exp|B_x\cap V_x)\i)*\ph^x \in \Om^p_{G_x-\text{hor}}(S_x)^{G_x}$$
and by lemma \nmb!{4.2} this form corresponds uniquely to a
differential form
$\om^x\in \Om^p_{\text{hor}}(G.S_x)^G$
which satisfies $(i|\Si_x)^*\om^x = \om|\Si_x$,
since the exponential mapping commutes with the respective
restriction mappings.
Now the intersection $G.S_x\cap \Si$ is the disjoint union of all the
open sets $w_j(\Si_x)$ where we pick one $w_j$ in each left coset of
the subgroup $W(\Si)_x$ in $W(\Si)$. If we choose $g_j\in N_G(\Si)$
projecting on $w_j$ for all $j$, then
$$\align
(i|w_j(\Si_x))^*\om^x &= (\ell_{g_j}\o i|\Si_x \o w_j\i)^*\om^x \\
&= (w_j\i)^*(i|\Si_x)^*\ell_{g_j}^*\om^x \\
&= (w_j\i)^*(i|\Si_x)^*\om^x = (w_j\i)^*(\om|\Si_x) = \om|w_j(\Si_x),
\endalign$$
so that $(i|G.S_x\cap \Si)^*\om^x = \om|G.S_x\cap \Si$.
We can do this for each point $x\in \Si$.

Using the method of Palais (\cit!{16}, proof of 4.3.1\ign{ 5.8, 5.10})
we may find a sequence of
points $(x_n)_{n\in \Bbb N}$ in $\Si$ such that the $\pi(\Si_{x_n})$
form a locally finite open cover of the orbit space
$M/G\cong \Si/W(\Si)$, and a smooth partition of unity $f_n$
consisting of $G$-invariant functions with
$\operatorname{supp}(f_n)\subset G.S_{x_n}$.
Then $\tilde\om:= \sum_n f_n \om^{x_n}\in\Om^p_{\text{hor}}(M)^G$ has
the required property $i^*\tilde\om=\om$. \qed

\head\totoc\nmb0{5}. Basic versus equivariant cohomology \endhead

\subhead\nmb.{5.1}. Basic cohomology \endsubhead
For a Lie group $G$ and a smooth $G$-manifold $M$, by \nmb!{2.2} we
may consider the basic cohomology
$H_{G-\text{basic}}^p(M)=H^p(\Om^*_{\text{hor}}(M)^G,d)$.

\subhead\nmb.{5.2}. Equivariant cohomology, Borel model \endsubhead
For a topological group and a topological $G$-space the equivariant
cohomology was defined as follows, see \cit!{3}:
Let $EG\to BG$ be the classifying $G$-bundle, and consider the
associated bundle $EG\x_GM$ with standard fiber the $G$-space $M$.
Then the equivariant cohomology is given by $H^p(EG\x_GM;\Bbb R)$.

\subhead\nmb.{5.3}. Equivariant cohomology, Cartan model \endsubhead
For a Lie group $G$ and a smooth $G$-manifold $M$ we consider the
space
$$(S^k\g^*\otimes \Om^p(M))^G$$
of all homogeneous polynomial
mappings $\al:\g\to \Om^p(M)$ of degree $k$ from the Lie algebra $\g$
of $G$ to the space of $k$-forms, which are $G$-equivariant:
$\al(\Ad(g\i)X)=\ell_g^*\al(X)$ for all $g\in G$.
The mapping
$$\gather
d_\g:A_G^q(M)\to A_G^{q+1}(M)\\
A_G^q(M):=\bigoplus_{2k+p=q}(S^k\g^*\otimes \Om^p(M))^G\\
(d_\g\al)(X):= d(\al(X))-i_{\ze_X}\al(X)
\endgather$$
satisfies $d_\g\o d_\g=0$ and the following result holds.

\proclaim{Theorem} Let $G$ be a compact connected Lie
group and let $M$ be a smooth $G$-manifold. Then
$$H^p(EG\x_GM;\Bbb R) = H^p(A^*_G(M),d_\g).$$
\endproclaim

This result is stated in \cit!{1} together with some arguments, and
it is attributed to \cit!{5}, \cit!{6} in chapter 7 of \cit!{2}.
I was unable to find a satisfactory published proof.

\subhead\nmb.{5.4}. \endsubhead
Let $M$ be a smooth $G$-manifold. Then the obvious embedding
$j(\om)=1\otimes \om$ gives a mapping of graded differential algebras
$$j:\Om^p_{\text{hor}}(M)^G\to (S^0\g^*\otimes \Om^p(M))^G
\to \bigoplus_{k}(S^k\g^*\otimes \Om^{p-2k}(M))^G = A^p_G(M).$$ On
the other hand evaluation at $0\in\g$ defines a homomorphism of
graded differential algebras
$\operatorname{ev}_0:A^*_G(M)\to \Om^*(M)^G$, and
$\operatorname{ev}_0\o j$ is the embedding $\Om^*_{\text{hor}}(M)^G
\to \Om^*(M)^G$. Thus we get canonical homomorphisms in cohomology
$$\CD
H^p(\Om^*_{\text{hor}}(M)^G) @>{J^*}>> H^p(A^*_G(M),d_\g) @>>>
     H^p(\Om^*(M)^G,d)\\
@|  @|  @| \\
H^p_{G-\text{basic}}(M) @>>> H^p_G(M) @>>> H^p(M)^G.
\endCD$$
If $G$ is compact and connected we have $H^p(M)^G=H^p(M)$, by
integration and homotopy invariance.

\enddocument